\documentclass[a4paper,11pt]{article}

\RequirePackage{fix-cm}

\usepackage{amsmath,amsfonts,amssymb}
\usepackage{a4wide}
\usepackage{graphics}
\usepackage{epsfig}
\usepackage{epic,eepic}
\usepackage{lscape}
\allowdisplaybreaks

\newcommand{\R}{{\mathbb{R}}}
\newcommand{\C}{{\mathbb{C}}} 
\newcommand{\N}{{\mathbb{N}}} 
\newcommand{\ds}{\displaystyle}
\newcommand{\la}{\langle}
\newcommand{\ra}{\rangle}

\renewcommand{\Im}{\mathop{\rm Im}\nolimits}
\renewcommand{\Re}{\mathop{\rm Re}\nolimits}

\newenvironment{proof}{\noindent\textbf{Proof.}  }{\hspace*{\fill}$\Box$ \\[2mm]}

\newtheorem{proposition}{Proposition}[section]
\newtheorem{theorem}[proposition]{Theorem}
\newtheorem{corollary}[proposition]{Corollary}

\newtheorem{example}[proposition]{Example}

\title{On the Use of Elliptic Regularity Theory for the Numerical Solution
of Variational Problems}

\author{Axel Dreves  \& Joachim Gwinner \& Nina Ovcharova 
\thanks{Universit\"at der Bundeswehr M\"unchen,
Department of Aerospace Engineering,
Werner-Heisenberg-Weg 39,
85577 Neubiberg/Munich, Germany, 
Axel.Dreves@unibw.de, Joachim.Gwinner@unibw.de., Nina.Ovcharova@unibw.de}}

\date{\today}

\begin{document}

\maketitle

\begin{abstract}
In this article we show the crucial role of elliptic regularity theory
for the development of efficient numerical methods for the solution of
some variational problems.
Here we focus to a class of elliptic multiobjective optimal control problems 
that can be formulated as jointly convex generalized Nash equilibrium problems (GNEPs)
and to nonsmooth boundary value problems that stem from contact mechanics leading to elliptic variational inequalities (VIs).

\end{abstract}

\noindent {\bf Keywords:} 
complementarity problem, dual mixed formulation, elliptic boundary value problem,
jointly convex generalized Nash equilibrium problem, Lagrange multiplier,
multiobjective optimal control, normalized Nash equilibrium, obstacle problem, 
saddle point formulation, Signorini problem, smooth domain, unilateral contact,
variational inequality.

\section{Introduction}
As noted in the survey paper \cite{MR1251844}, elliptic regularity theory is of essential importance for the derivation of error estimates of the finite element method (FEM) for the numerical solution of nonsmooth boundary value problems formulated as variational inequalities. This is now well documented in the literature starting from the pioneering work of Falk \cite{MR0391502}. More recent examples of this research direction are the 
paper \cite{MR2982470} on the $h$-FEM treatment of unilateral crack problems and other nonsmooth constraints and the paper \cite{MR3061074} on $hp$-FEM convergence for unilateral contact problems with Tresca friction in plane linear elastostatics.

This article is concerned with other applications of elliptic regularity theory.
First we consider a class of elliptic multiobjective optimal control problems 
formulated as jointly convex generalized Nash equilibrium problems.
As will be detailed below, a rather straightforward variational formulation of such a problem
leads to a generalized Nash equilibrium problem (GNEP), where however each player has to satisfy different constraints that depend on the control of the other players. Thus one obtains more involved quasi-variational inequalities, in contrast to variational inequalities that can be obtained 
when considering normalized solutions of jointly convex GNEPs as is shown in the recent paper \cite{ADJG-JOTA}, based on regularity of the underlying elliptic boundary value problem.

Then we turn to Signorini mixed boundary value problems, unilateral frictionless contact problems and other nonsmooth boundary value problems that can be formulated as variational inequalities with a coercive bilinear form. To get rid of relatively complicated constraints as e.g., inequality constraints and to obtain simpler nonnegativity constraints or box constraints one can introduce Lagrange multipliers similar as in constrained optimization in finite dimensions. In addition to simplification for better numerical treatment, there is also an intrinsic interest in Lagrange multipliers as dual variables. Often in applications they have a clear physical meaning and are more of interest than the primal variables; speaking in the language of continuum mechanics, the engineer is often more interested in the stresses and strains than in the displacements. This motivates  multifield variational formulations and multiple saddle point problem formulations, see 
\cite{MR1975268,MR3049426,hammamet}. While for linear elliptic boundary value problems the passage from the primal variational formulation to a dual mixed formulation or a saddle point problem form involving a Lagrange multiplier is a standard procedure and while there are the well-established mixed finite element methods \cite{MR1115205,MR1958060} for their numerical treatment, such a procedure for non-smoothly constrained problems has to overcome several difficulties. 
First, the standard approach to existence of Lagrange multipliers for inequality constrained optimization in infinite dimensional spaces relies on the Hahn-Banach separation theorem and needs an interior point condition (Slater condition) with respect to the ordering cone in the image space. However, the topological interior of such an ordering cone in standard function (Hilbert or Banach) spaces, as e.g. the interior of the cone 
$L^p_+$ of non-negative $L^p$ functions is empty. So one may resort to the nonempty quasi relative interior of $L^p_+$ and one may impose a Slater-like condition, that is, the existence of a feasible point that lies in the quasi relative interior of $L^p_+$. However, as a counterexample of Daniele and Giuffr{\`e} \cite{MR2340691}
shows, this condition is not sufficient, and extra more complicated assumptions 
or related involved conditions that are actually equivalent  are needed to ensure the existence of a Lagrange multiplier, 
see \cite{MR2438595,MR2340691,MR2317768,MR3245967}.

Therefore we proceed in another way and show how by a simple formula one obtains a Lagrange multiplier in the dual of the preimage space thus even reducing the variational inequality to a complementarity problem. By this simple approach,    
the Lagrange multiplier lives in the dual of the Sobolev space of the variational problem, thus at first, is a general measure which may be singular.  Here regularity theory comes into play to conclude that the Lagrange multiplier is indeed an $L^p$ function. Thus from an inequality constraint, one finally obtains a Lagrange multiplier in the cone $L^p_+$ of non-negative $L^p$ functions. This approach works also with not necessarily symmetric bilinear forms, when the equivalence to convex quadratic optimization is lost; it even works for  nonlinear operators. Moreover, we can combine such dual mixed formulations for variational problems with inequality constraints  via non-negative Lagrange multipliers with mixed formulations for variational inequalities of second kind where the Lagrange multiplier is in a simple box set. This applies to unilateral contact problems with Tresca friction.
 
The outline of this article is as follows. The next section provides a review of elliptic regularity theory dealing with the linear Dirichlet problem, the scalar unilateral boundary value problem (obstacle problem), and  frictionless unilateral contact contact of linear elastostatics. In section 3 we consider a class of elliptic multiobjective optimal control problems and show following \cite{ADJG-JOTA}
how based on elliptic regularity theory, these problems can be reformulated 
as so-called  jointly convex generalized Nash equilibria.
In  section 4 we present a direct approach to mixed formulations of some nonsmooth variational problems and of associated variational inequalities. 
The article ends with some conclusions and an outlook to some open problems.

\section{A Review of Elliptic Regularity Theory}
In this section we review the elliptic regularity theory that is needed for the understanding of the subsequent sections. 

\subsection{Regularity of linear scalar Dirichlet problem}
In this subsection we are concerned with the regularity of the solution of the Dirichlet problem with a (scalar) linear second-order elliptic operator $L$; that is, the operator $L$
is of the form (summation convention employed) 
$$
Lu = D_i(-a^{ij}(x)~D_j u) +  a(x) u \,,
$$
where the coefficients $a^{ij}~(i,j=1, \ldots,d)$ and $a$ are assumed to be 
bounded, measurable functions on a domain $\Omega \subset \R^d$ and moreover,
$a$ is non-negative and  there exists a positive number $\alpha$ such that 
$$
a^{ij}(x)~\xi_i\xi_j \geq \alpha |\xi|^2 \,, \forall x \in \Omega, \xi \in \R^d \,.
$$
A simple example is $L= - \Delta$, the negative Laplacian on $\R^d$; on the other hand, 
lower order terms involving $D_i u$ can easily be included 
in the definition of $L$. The operator $L$ above gives rise to the bilinear form
$$
{\cal L}(u,v) = \int_{\Omega} [a^{ij}(x) D_j u D_i v + a(x) uv]~ dx \,.
$$
Let in addition $f$ be (locally) integrable on $\Omega$ and $\varphi$
belong to $H^1(\Omega)$, the Sobolev space of all $L^2$ functions on $\Omega$
with weak $L^2$ derivatives, see \cite{MR2424078}.
Then a function $u \in H^1(\Omega)$ is called a weak solution of the Dirichlet problem:
$$
Lu = f, \quad u= \varphi \mbox{ on } \partial \Omega,
$$
if $u - \varphi \in H_0^1(\Omega)$ and $u$ satisfies
$$
{\cal L}(u,v) = \int_{\Omega} f v ~dx , \, \forall v \in C_0^1(\Omega) \,.
$$

The following example of a domain with a reentrant corner taken from the book of Braess
\cite{MR2322235}
shows that even for smooth data $f, \varphi$ we cannot expect the solution to be 
in $H^2 (\Omega)$, not to mention in  $C^2 (\Omega)$, what is suggested by a classic treatment of partial differential equations.

\begin{example}
Let 
$$
\Omega =\{ x = (x_1,x_2) \in \R^2~:~ |x| < 1, x_1 < 0 \mbox{ or } x_2 > 0 \} \,.
$$
Identify $\R^2$ with $\C$. Let $z = x_1 + i x_2 = \rho \exp(i \theta)$ and consider
$$
w(z) = z^{2/3}; \quad u(x) = \Im w(z)   =  \rho^{2/3} \sin \left (\frac{2}{3}\theta \right ) \,.
$$
So $u$ is harmonic and  $u \in H^1(\Omega)$ solves 
\begin{eqnarray*}
\begin{array}{rccl}
 \Delta u &=& 0 \qquad \:\: & \mbox{ in } \Omega; \\
    u(\exp(i \theta)) & = & \sin \left (\frac{2}{3}\theta \right )
		  & \mbox{ for } 0 \leq \theta \leq \frac{3}{2} \pi , \\
   u &=& 0  \qquad \:\: & \mbox{ elsewhere on } \partial \Omega  \,.
	\end{array} 
\end{eqnarray*}
Since $ w'(z) = \frac{2}{3}\ z^{-\frac{1}{3}}$, even the first derivatives of $u$ are not 
bounded for $z \rightarrow 0$.
\end{example}

There are two options for a domain to obtain regularity $u \in H^2 (\Omega)$:
smoothness of the boundary $\partial \Omega$ or convexity of the domain. For the first let us recall from the monograph of Gilbarg and Trudinger
\cite[Theorem 8.12]{MR1814364}

\begin{theorem}
Suppose that $\partial \Omega$ is of class $C^2$.  Moreover assume the coefficients 
$a^{ij}$ are uniformly Lipschitz continuous in $\Omega$ and for the data assume
$f \in L^2(\Omega)$ and $\varphi \in H^2(\Omega)$ such that 
$u - \varphi  \in H_0^1(\Omega)$ with a weak solution $u$ of the above Dirichlet problem.
Then also $u \in H^2(\Omega)$.
\end{theorem}
For such a regularity result and for its direct proof we can also refer to the monograph of Aubin \cite[chapter 7, section 1-8, Theorem 1-1]{MR0478662} and to the monograph of Kinderlehrer and Stampacchia \cite[chapter IV, appendix A]{MR1786735}.

Regularity results for the Dirichlet problem for elliptic operators, respectively for the Laplacian  on convex domains and on more general so-called semiconvex domains (here a bounded domain is semiconvex, if for any $x \in \partial \Omega$ there exists an open ball $B_x \subset \R^d \setminus \bar{\Omega}$
 with $\bar{B_x} \cap \bar{\Omega} = \{x\} $) 
are established in the work of Kadlec \cite{MR0170088}  and of Mitrea et al.\,\cite{MR2593333}. Let us also mention the regularity results for solutions of the equations of linear elasticity in convex plane polygonal domains by Bacuta and Bramble
\cite{MR2019187}.

\subsection{Regularity of the scalar unilateral boundary value problem}
Let us turn to the regularity of scalar unilateral boundary value problems, in particular of Signorini boundary value problems. We are also concerned with the regularity of domain obstacle problems, since domain obstacle and boundary obstacle (Signorini) problems are related as follows.

Let $\Gamma_D, \Gamma_S$ be two disjoint smooth and open subset of $\partial \Omega$ such that
$\partial \Omega = \bar \Gamma_D \cup \bar \Gamma_S$. Let $A$ be a linear elliptic operator
defined by $Au = - D_j(a^{ij} D_i u)$ with coefficients $a^{ij}$ as above giving the bilinear form 
$$
a(u,v) = \int_\Omega a^{ij}(x)~ D_i u (x) ~ D_j v(x) ~dx \,.
$$
Let $\psi \in H^1(\Omega)$ with $\psi \leq 0$ on $\Gamma_S$, let $\tilde \psi$ be the unique solution of the Dirichlet problem
$$
\Delta \tilde \psi = f \mbox{ in } \Omega, 
\qquad \tilde \psi = \psi \mbox{ on } \partial \Omega
$$
and assume that $\tilde \psi \in H^2(\Omega)$.  
Let 
$$
V = \{ v \in H^1(\Omega): ~ v = 0 \mbox{ on } \Gamma_D \}
$$
and define the closed convex subsets of $V$:
\begin{eqnarray*}
K &=&  \{ v \in V: ~ v \geq \psi \mbox { on } \Gamma_S \}, \\
\tilde K &=& \{ v \in V: ~ v \geq \tilde \psi \mbox{ in } \Omega,
 v = \psi \mbox { on } \Gamma_S \} \,.
\end{eqnarray*}
Then there holds the following

\begin{theorem} \label{equivVI}
If $u$ is the solution of the VI (domain obstacle problem)
$$
u \in \tilde K, \quad  a(u,v-u) \geq \int_\Omega f(v-u) ~ dx \quad \forall v \in \tilde K,
$$
then $u$ resolves the VI (Signorini problem)
$$
u \in  K, \quad a(u,v-u) \geq \int_\Omega f(v-u) ~ dx \quad \forall v \in  K \,.
$$
\end{theorem}
For its proof see the proof of Theorem 9.3 in \cite[chapter IV]{MR1786735}.

In virtue of Theorem \ref{equivVI} we can conclude from the regularity result
\cite[chapter IV, Theorem 2.3]{MR1786735} for the domain obstacle problem the
following regularity result for the Signorini problem with $A = - \Delta$:

\begin{theorem} \label{regulSig}
Suppose $f \in L^s(\Omega)$ and $\max(- \Delta \tilde \psi -f,0) \in L^s(\Omega)$
for some $s > d$. Then the solution of the above Signorini problem with
$A = - \Delta$ lies in $H^{2,s}(\Omega) \cap C^{1,\lambda}(\bar \Omega),
\lambda = 1 - (d/s)$. Hence $\Delta u \in L^s(\Omega)$.
\end{theorem}

There is a refinement concerning the regularity of the domain obstacle problem at the boundary by Jensen \cite{MR578201}. He has proven the local regularity 
$W^{2,\infty}$ of the solution at boundary points. However, Kinderlehrer
\cite{MR656002} has provided the following example of a scalar Signorini problem
with a solution that fails to be in $H^2(\Omega)$. Here even the boundary obstacle is zero,
but $\partial \Gamma_D \cap \partial \Gamma_S \not= \emptyset$ with the 
Dirichlet part $\Gamma_D$ and the Signorini part $\Gamma_S$.

\begin{example}
Let 
$$
\Omega =\{ x = (x_1,x_2) \in \R^2~|~ |x| < 1, x_2 > 0 \}
$$
with the mutually disjoint, open boundary parts 
\begin{eqnarray*}
 && \Gamma_S = \{ x = (x_1,0) ~|~ -1 < x_1 < 0 \}, \\
 && \Gamma_N = \{ x = (x_1,0) ~|~ 0 < x_1 < 1 \}, \\
 && \Gamma_D = \{ x = (x_1,x_2) ~|~ |x| =1, x_2 > 0 \} \,.
\end{eqnarray*}
Let $z = x_1 + i x_2 = \rho \exp(i \theta) $ and consider
$$
u(x) = - \Re z^{1/2}  = - \rho^{1/2} \cos (\theta/2)  \,.
$$
So $u$ is harmonic and $u(x) = 0$ for $x_1 < 0, x_2 = 0$. By the Cauchy-Riemann
differential equations,
\[ \ds \frac{\partial}{\partial \nu}~u(x_1,0) = - \ds \frac{\partial}{\partial x_2}~u(x_1,0) = - \ds \frac{\partial}{\partial x_1}~\Im z^{1/2} =
\left\{ \begin{array}{cc}
0 & \mbox{if } x_1 > 0 \,,\\ \ds \frac{1}{2} |x_1|^{-\frac{1}{2}} & \mbox{if } x_1 < 0 \,,
\end{array} \right.
\, x_2 = 0 \,. \]
Hence $u \in H^1(\Omega)$ satisfies $- \Delta u = 0$ in $\Omega$ and the Neumann, respectively Dirichlet boundary conditions
$$
\ds \frac{\partial u}{\partial \nu} = 0 \mbox{ on } \Gamma_N, \,
 u = - \ds \cos \frac{\theta}{2} \mbox{ on } \Gamma_D \,,
$$
and the Signorini boundary conditions
$$
u ~\ds \frac{\partial u}{\partial \nu} = 0, \, u \geq 0, \, 
\ds \frac{\partial u}{\partial \nu} \geq 0 \mbox{ on } \Gamma_S \,.
$$
Thus $u$ solves the VI
$$
u \in K, \quad \int_\Omega \nabla u \cdot \nabla(v - u)~dx \geq 0, \, \forall v \in K \,,
$$
where
$$
K = \{ v \in H^1(\Omega)~|~ v \geq 0 \mbox{ on } \Gamma_S,\, 
v = - \ds \cos \frac{\theta}{2} \mbox{ on } \Gamma_D \} \,.
$$
Note that $\ds \frac{\partial^2}{\partial x_1^2} 
= \cos^2 \theta ~ \ds \frac{\partial^2}{\partial \rho^2} + \ldots$,  
$\ds \frac{\partial^2}{\partial x_2^2} 
= \sin^2 \theta ~ \ds \frac{\partial^2}{\partial \rho^2} + \ldots$,
$$
\int \!\!\! \int_{\Omega} |u_{\rho \rho}|^2 ~dx =
(1/16)~ \int_0^\pi \! \int_0^1 \rho^{-3}~\rho~  \cos^2(\theta/2) d\rho \, d\theta;
\int_0^1\rho^{-2}~d\rho = \infty \,
$$
so $u$ cannot lie in $H^2(\Omega)$.
\end{example}

\subsection{Variational formulation of frictionless unilateral 
contact problem of linear elastostatics}

Before we continue our review of elliptic regularity theory
addressing frictionless unilateral contact problems we introduce some notation from continuum mechanics and describe the variational form
of unilateral contact problems as variational inequalities (of first kind, following the terminology of \cite{MR2423313}).

Let us assume Hooke's  law and small deformations of a non--homogeneous, anisotropic body. For notational simplicity we focus to the case of 
plane elasticity; the three-dimensional case poses no additional difficulty in deriving the variational formulation. So let $ \Omega \subset \R^2$ be a bounded plane domain with Lipschitz boundary  $\Gamma$ ($\Gamma \in C^{0, 1}$),  occupied by an
  elastic body,  and let 
 $ \underline{x} = (x_1 , x_2) $ be a Cartesian coordinate system. Then
 $\underline{n} = (n_1 , n_2)$, the unit outward normal to $\Gamma$, exists almost  everywhere and $ \underline{n} \in [L^\infty (\Gamma)]^2$, see e.g. 
 \cite[Theorem 5.4]{MR961258}. \\ With the displacement vector 
$ \underline{v} =  (v_1 , v_2)$ to lie in the Sobolev space 
$ [H^1 (\Omega)]^2 $ the  linearized strains are given by
 \begin{equation}\label{E3.1} \varepsilon_{ij} (v) = \frac{1}{2}
 \left( \frac{\partial v_i}{\partial x_j} + \frac{\partial v_j}{\partial x_i}
 \right) \quad ( i, j = 1, 2) \end{equation}
 and Hooke's law relating strains and stresses reads
 \begin{equation}\label{E3.2} \tau_{ij} = E_{ij k l} \:\: \varepsilon_{k l}
 \quad (i, j = 1, 2) \: , \end{equation}
 where we use the summation convention over a repeated index within the range
 $1 , 2$ and where the elasticity coefficients 
$ E_{ij k l} \in L^\infty (\Omega)$ satisfy
 \[ E_{ij k l}  = E_{k l ij}  = E_{ji k l} \: ; \]
 \begin{equation}\label{E3.3} \exists c_0 > 0: \: E_{ij k l} \:\:
 \varepsilon_{ij} \:\: \varepsilon_{k l} \: \geq \: c_0 \:
 \varepsilon_{ij} \:\: \varepsilon_{ij} \quad \forall \: \varepsilon_{ij}
 = \varepsilon_{ji} \: . \end{equation}
 With the given vector 
$ \underline{F} = (F_1 , F_2) \in [L^2 (\Omega)]^2 $ the stress
 field has to satisfy the equilibrium equations
 \begin{equation}\label{E3.4} \frac{\partial \tau_{ij}}{\partial x_j}  + F_i
 = 0 \quad (i = 1, 2) \: . \end{equation}
 The traction vector $\underline{b}$ on the boundary, where
 \[ b_i = \tau_{ij} \: n_j \]
 can be decomposed into the normal component
 \[ b_n = b_i \: n_i = \tau_{ij} \: n_i \: n_i \]
 and the tangential component
 \[ b_t = b_i \: t_i = \tau_{ij} \: t_i \: n_j \: , \]
 where $ \underline{t} = (t_1 , t_2) = (- n_2 , n_1 ) $ is the unit tangential
 vector. Likewise the displacement $v$ can be decomposed (see 
 \cite[Chapter 5]{MR961258}, \cite{MR1301021} for the relevant trace theorems):
 \[ v_n = v_i n_i \: , \:\: v_t = v_i  t_i \: . \]
 To describe the boundary conditions, let 
$ \Gamma = \overline{\Gamma}_{D} \cup
 \overline{\Gamma}_{N} \cup \overline{\Gamma}_{S} $, where the open parts $ \Gamma_D , \Gamma_N , $ and $ \Gamma_S$ are mutually disjoint. 
Eventually nonzero displacements
$ \underline{D} \in [ H_1 ( \Gamma_D) ]^2$,
 respectively tractions $ \underline{T} \in [ L_2 ( \Gamma_N) ]^2$ are prescribed on  $\Gamma_D$, resp. $ \Gamma_N$, i.e.,
 \begin{eqnarray}
\label{E3.5} v_i & = & D_i \qquad  \mbox{on} \quad  \Gamma_D \: , \\
\label{E3.6} b_i & = & T_i  \qquad \: \mbox{on} \quad  \Gamma_N \: ,
\end{eqnarray}
whereas on $\Gamma_S$ the frictionless unilateral contact conditions
(Signorini's conditions for $v_n$ and $b_n$)
\begin{equation}\label{E3.7} v_n \leq g , \:  b_n \leq 0 , \: 
(v_n -g)\, b_n = 0, \:
b_t = 0 \end{equation}
with a given gap function $g \in L_2(\Gamma_S)$ 
are imposed. To make the contact problem
meaningful we assume  $ meas (\Gamma_S) > 0$. Here we also require 
$ meas(\Gamma_D) > 0 $, hence rigid body motions are excluded and the variational problem becomes coercive. \\
Now the problem (\ref{E3.1}), (\ref{E3.2}), (\ref{E3.4}) - (\ref{E3.7}) 
can be formulated as the following variational inequality (VI): Find
$u \in K $ such that
\begin{equation}\label{E2.1}
\beta (u, v-u) \geq \lambda (v-u) \quad \forall v \in K \: .
\end{equation}
where we introduce the bilinear form, respectively the linear form
\begin{eqnarray*} \beta (\underline{v}, \underline{w})
 & =  & \int_{\Omega} E_{i j k l} \:
\varepsilon_{ij}  (\underline{v}) \:
\varepsilon_{k l}  (\underline{w}) \: dx \: , \\
\lambda (\underline{v}) & = & \int_{\Omega} F_i v_i \: dx + \int_{\Gamma_{N}} T_i v_i \:ds
\end{eqnarray*}
on the function space
\[ V = \left\{ \underline{v} \in [H^1 (\Omega) ]^2 \:  | \:  \underline{v} = 0 \mbox{ on } \Gamma_D \right\} \]
and the convex closed subset 
\[ K = \left\{ \underline{v} \in V \: | \: 
v_i = D_i \mbox{ on } \Gamma_D; \,
v_n \leq g \mbox{ on } \Gamma_S \right\} \: . \]

One may reduce the inhomogeneous inequality constraint $v_n = v_i n_i \leq g$ to the homogeneous inequality constraint 
$\tilde v_n = \tilde v_i n_i \leq 0$, thus simplify to a convex cone constraint
by subtraction of some appropriate extension $\underline{g}$ of $g \in L_2(\Gamma_S)$
to $[H^1 (\Omega) ]^2$.  However, this simple reduction for unilateral
constraints  does not work with more general bilateral constraints of the form 
$g_a \leq v_n = v_i n_i \leq g_b $, when the extended real-valued 
boundary obstacles have domains that intersect, i.e.,
$\mbox{dom } g_a \cap \mbox{dom } g_b \not= \emptyset$,
in particular in a three-dimensional situation.  

\subsection{Regularity of frictionless unilateral 
contact problem of linear elastostatics}

In view of his example given above Kinderlehrer \cite{MR656002}
could prove by a difference quotient technique 
that the solution $u$ of the Signorini problem is in $H^2$ except perhaps near points of
$ \partial \Gamma_S \cup \partial \Gamma_N \cup \partial \Gamma_D 
\subset \Omega$, more precisely the following result
for the $d-$dimensional mixed Signorini boundary value problem
in the case $g = 0$, what is  by the remark above, no loss of generality
concerning regularity.

\begin{theorem} \cite[Theorem 2.2]{MR656002}
Suppose for the data  
$ \underline{F}  \in [L^2 (\Omega)]^d $, 
$ \underline{T} \in [ H^1 ( \Gamma_D) ]^d$,
$ \underline{D} \in [ H^2 ( \Gamma_D) ]^d$.
Set $\Omega_\delta = \{ x \in \Omega~|~
\mbox{dist } (x, \partial \Gamma_S \cup \partial \Gamma_N
\cup \partial \Gamma_D) > \delta \}$ for $\delta > 0$.
Then for each  $\delta > 0$ there hold $ u \in (H^2(\Omega_\delta))^d$ 
and the Signorini conditions (\ref{E3.7}) pointwise a.e. on  $ \Gamma_S$.  
\end{theorem}

Sobolev imbedding of $H^2$ in spaces of H\"older continuous functions
implies the 

\begin{corollary}
Under the assumptions of the data as in the above theorem,
there holds
\begin{eqnarray*}
 \mbox{for }& d =2, & u\in [C^{0,\alpha} (\bar \Omega_\delta)]^2
 \mbox{ for  some } 0 < \alpha < 1  \,, \\
\mbox{and for }&  d = 3, & u\in [C^{0,\frac{1}{2}} (\bar \Omega_\delta)]^2 \,.
\end{eqnarray*}
\end{corollary}

By the theory of pseudodifferential operators Schumann \cite{MR986184},
extended the latter result to  $[C^{1,\alpha} (\Omega \cup \Gamma)]^2$
regularity of the solution $u$; the precise value of $\alpha$ is not  known.
To conclude this section, we refer to the survey \cite{MR1864872} of schumann who gives an excellent overview of the mathematical methods to prove regularity results for variational inequalities and unilateral problems in elasticity.
  
\section{From elliptic multiobjective optimal control 
to jointly convex generalized Nash equilibria}

In this section we consider a class of elliptic multiobjective optimal control problems and show following \cite{ADJG-JOTA}
how based on elliptic regularity theory, these problems can be reformulated 
as so-called  jointly convex generalized Nash equilibria.

\subsection*{The concept of jointly convex generalized Nash equilibrium problems} 

Let $V_\nu, \nu=1,\ldots,N$ be real separable Hilbert spaces or more general
reflexive, separable Banach spaces endowed with norms $\|\cdot \|_\nu,$ and define $V := V_1 \times \ldots \times V_N.$ Further, let $X$ be a nonempty, closed, and convex subset of $V$ and assume that the objective functions $\theta_\nu:V_1 \times \ldots \times V_N \to \mathbb R,$ 
$\theta_\nu(\cdot,x^{-\nu}) : V_\nu \to \mathbb R$ are convex for any fixed $x^{-\nu},$ where we use the notation $x=(x^1,\ldots,x^N)=(x^\nu,x^{-\nu})$ to emphasize the role of the variable $x^\nu$, but this notation does not mean a permutation.
In this setting the infinite dimensional jointly convex generalized Nash equilibrium (GNEP for short) has the following form
\begin{equation} \label{GNEP}
\min_{x^\nu} \theta_\nu(x^\nu,x^{-\nu}) 
\quad \text{subject to (s.t.)} \quad (x^\nu,x^{-\nu}) \in X 
\end{equation}
for all $\nu=1,\ldots,N.$ The reason for calling this problem jointly convex is that the strategies must belong to a common convex set $X,$ instead of each player having his own strategy set $X_\nu(x^{-\nu})$ depending on the rivals' strategy $x^{-\nu}.$
We call $\bar x$ a {\it generalized Nash equilibrium}, if $\bar x \in X$ satisfies 
\[
\theta_\nu(\bar x^{\nu},\bar x^{-\nu}) \le \theta_\nu(x^\nu,\bar x^{-\nu})\,, \quad \forall (x^\nu,\bar x^{-\nu}) \in X 
\]
for all $\nu=1,\ldots,N.$ Note that the concept of GNEPs goes back to the 
$1954$ paper \cite{MR0077069} of Arrow and Debreu.

Next let us introduce the Nikaido-Isoda function 
\[
\Psi(x,y) := \sum_{\nu=1}^N \left [ \theta_\nu(x^\nu,x^{-\nu}) - \theta_\nu(y^\nu,x^{-\nu}) \right ],
\]
see the $1955$ paper \cite{MR0073910} of Nikaido and Isoda, to define normalized solutions of a jointly convex GNEP:
$\bar x \in X$ is called a {\it normalized Nash equilibrium}, or a 
{\it normalized solution of the jointly convex GNEP} if 
\[
 \Psi(\bar x,y) \leq 0 \,, \quad \forall y \in X \,.
\]
Thus we get a characterization of some solutions, namely the normalized solutions of jointly convex GNEPs via a variational inequality in contrast to more involved quasi-variational inequalities that characterize the solutions of GNEPs in general form, not necessarily jointly convex. Therefore, computing normalized solutions of jointly-convex GNEPs is typically much easier than obtaining solutions of GNEPs in general form.
Since for a normalized Nash equilibrium $\bar x$  we have for all $\nu=1,\ldots,N$ and all
$(y^\nu,\bar x^{-\nu}) \in X$
\begin{eqnarray*}
\theta_\nu(\bar x^\nu,\bar x^{-\nu}) - \theta_\nu(y^\nu,\bar x^{-\nu}) =
\Psi(\bar x,(y^\nu,\bar x^{-\nu})) \leq \sup_{y \in X} \Psi(\bar x,y) \le 0,
\end{eqnarray*}
every normalized solution is also a generalized Nash equilibrium, i.e., 
for all $\nu=1,\ldots,N$ it holds that
\[
\theta_\nu(\bar x^\nu,\bar x^{-\nu}) \le \theta_\nu(y^\nu,\bar x^{-\nu}) \qquad \forall (y^\nu,\bar x^{-\nu}) \in X \,; 
\]
however, the converse is not true. 

\subsection*{Primal formulation of elliptic multiobjective optimal 
control problems} 

Let $\Omega \subset  \mathbb{R}^{d} \,(d=2,3) $ be a bounded Lipschitz domain and let $V := H_0^1(\Omega)$ denote the Sobolev space of all 
$L^2$ functions on $\Omega$ with weak $L^2$ derivatives and zero boundary values. 
Let $U^\nu := L^2(\Omega)$ be the space for the controls $u^\nu$ for all $\nu=1,\ldots,N.$  
We have the weights $\gamma_\nu >  0, \beta_\nu >0$, the given data  $f, g^\nu \in L^2(\Omega), \, (\nu=1,\ldots,N); a^\nu,b^\nu \in \mathbb R$ with $a^\nu \le b^\nu$ 
($\nu=1,\ldots,N$), $a^0,b^0 \in H^1(\Omega)$ with $a^0(x) < b^0(x)$ and some continuous, compact, and linear operators $\chi^\nu:V \to L^2(\Omega)$ ($\nu=1,\ldots,N$). 
Then we consider the following problem
\begin{eqnarray*} 
(I) & \min\limits_{y,u^\nu} &\frac{1}{2} \left \| \chi^\nu y - g^\nu \right \|^2_{L^2(\Omega)} + \frac{\gamma_\nu}{2} \left \|u^\nu \right \|^2_{L^2(\Omega)} 
\\ & \text{ s.t. } & 
L y = \sum_{\mu=1}^N \beta_\mu u^{\mu} + f, \quad y|_{\partial \Omega} = 0, \\
& & a^0(x) \le y(x) \le b^0(x), \text{ a.e. in } \Omega,\\
& & a^\nu \le  u^\nu(x) \le b^\nu, \qquad \, \text{ a.e. in } \Omega,
\end{eqnarray*}
for all $\nu=1,\ldots,N.$ In this problem every player $\nu$ minimizes his own cost function through his individual control variable $u^\nu$ and the common state variable $y$. The state is determined by the controls of all players via a partial differential equation 
(pde) given by a linear elliptic partial differential operator $L$ 
of second order as introduced in the previous section.

To provide a functional analytic meaning we can write the above pde constraint in variational form as 
\[
 y \in V: \quad 
{\cal L}(y,w) = 
\langle \sum_{\mu=1}^N \beta_\mu u^ \mu + f, w \rangle_{L^2(\Omega)} \,,
\quad \forall w \in V\,.
\]
Using the continuous embedding from $H_0^1(\Omega)$ in $L^2 (\Omega),$ the state constraints $a^0  \le y \le b^0$ are to be understood in the $L^2$ sense as the control constraints $a^\nu \le  u^\nu \le b^\nu, $ which imply that the controls are actually $L^\infty$ functions, since $a^\nu,b^\nu \in \mathbb R$. 

The elliptic multiobjective optimal control problem $(I)$ is, however, not a GNEP, since the state $y$ is a common optimization variable for all players. If we introduce different state variables $y^\nu$ for each player $\nu \in \{1,\ldots,N\},$ and if we could guarantee that all the states are equal, 
we get a GNEP. But we do not get a jointly convex GNEP,
since the players then have different constraints 
$L~ y^\nu = \sum_{\mu=1}^N \beta_\mu u^{\mu} + f,$ depending on the controls of the other players. For the numerical solution of these GNEPs in general form one can use its' optimality conditions that are equivalent to quasi-variational inequalities, and are much harder to solve than VIs. Also the number of algorithms for the solution of quasi-variational inequalities is rather limited.  
Therefore our next aim is to develop a jointly convex reformulation.

\subsection*{Reduced multicontrol formulation of the elliptic multiobjective optimal control problems} 

Since by the Lax-Milgram theorem $L: H_0^1(\Omega) \to H^{-1}(\Omega) $ is an isomorphism,
we can use the inverse $L^{-1}:H^{-1}(\Omega) \to H_0^1(\Omega)$ to define the multicontrol to state map
\[
S(u) := L^{-1} \left (\sum_{\mu=1}^N  \beta_\mu u^{\mu} + f \right ),
\]
and this is a continuous map affine linearly dependent on $(u^1,\ldots,u^N).$
 Since  $L^2(\Omega)$ is compactly embedded in $H^{-1}(\Omega),$ see \cite{MR2424078}, 
 this is even a completely continous map from $[L^2(\Omega)]^N$ to $H^1_0(\Omega).$
Hence, we obtain the equivalent reduced problem
\begin{eqnarray*} 
(II) & \min\limits_{u^\nu} & \frac{1}{2} \left \| \chi^\nu S(u^\nu,u^{-\nu}) - g^\nu \right \|^2_{L^2(\Omega)} + 
\frac{\gamma_\nu}{2} \left \|u^\nu \right \|^2_{L^2(\Omega)} \\ & \text{ s.t. }  &
 a^0(x) \le  S(u^\nu,u^{-\nu})(x) \le b^0(x), \text{ a.e. in } \Omega,\\
&& a^\nu \le u^\nu(x) \le b^\nu,  \qquad \qquad \quad \quad \, \text{ a.e. in } \Omega,
\end{eqnarray*}
for all $\nu=1,\ldots,N,$ which is a jointly convex GNEP. A similar problem was first considered in \cite{MR3098015} as a GNEP and using a penalty approach and a strict uniform feasible response assumption, the existence of a solution was shown. Further, using the Nikaido-Isoda function, this reformulation (II) was used in \cite{MR3397432} to show existence of a Nash equilibrium for the equivalent problem (I). Moreover it was shown that one can solve these reformulations (II) (even for parabolic and not only elliptic  pdes) via a primal-dual path-following method based on the Nikaiod-Isoda function.

Let us stress that (II) is already a jointly convex GNEP.
However, every evaluation of $S(u^\nu,u^{-\nu})$ requires the solution of the pde. To avoid this we give a third equivalent formulation to (I).

\subsection*{A multistate formulation of the elliptic multiobjective optimal 
control problems} 

Now we assume that $\partial \Omega$ is of class $C^2$ or $\Omega$ is convex. Then regularity theory for elliptic equations
with Dirichlet boundary conditions, as exposed in the previous section, 
guarantees that the solution $S(u)$ is even in $H^2(\Omega)$.
Therefore, 
\[
w^\nu :=  L^{-1} \left (\beta_\nu u^\nu \right )
\in H_0^1(\Omega) \cap H^2(\Omega) =:W
\]
for all $\nu=1,\ldots,N.$ These $w^\nu$ will become our new optimization variables. Indeed we now have 
\begin{equation} \label{eq:y}
y=S(u) = L^{-1} f + \sum_{\mu=1}^N w^\mu,
\end{equation}
and, since $L w^\nu \in L^2(\Omega)$, the equation
\begin{equation} \label{eq:u}
u^\nu = - \frac{1}{\beta_\nu} L w^\nu  
\end{equation}
holds in $L^2(\Omega)$ for all $\nu=1,\ldots,N.$ Thus we arrive at the equivalent problem
\begin{eqnarray*} 
(III) & \min\limits_{w^\nu \in W} & \frac{1}{2} \left \| \chi^\nu \left (\sum_{\mu=1}^{N} w^\mu \right ) + \chi^\nu  L^{-1} f- g^\nu \right \|^2_{L^2(\Omega)}
+ \frac{\gamma_\nu}{2 \beta_\nu^2} \left \| L w^\nu \right \|^2_{L^2(\Omega)} \\ & \text{ s.t. }  &
 a^0(x)  \le \left (L^{-1} f + \sum_{\mu=1}^N  w^\mu \right )(x) \le b^0(x),  \quad \text{ a.e. in } \Omega,\\ 
&& a^\nu \beta_\nu \le \left(L w^\nu \right )(x) \le b^\nu \beta_\nu, \qquad \qquad \text{ a.e. in } \Omega,
\end{eqnarray*}
for all $\nu=1,\ldots,N. $
Now, defining the common feasible set
\begin{eqnarray*} 
\widetilde {W}:= \Bigg \{ (w^1,\ldots,w^N) \in W^N &\Bigg |& a^\nu \beta_\nu \le \left(L w^\nu  \right )(x) \le b^\nu \beta_\nu \, \,
 (\forall \nu=1,\ldots,N), \\
&&  a^0(x)  \le \left (L^{-1} f + \sum_{\mu=1}^N  w^\mu \right )(x) \le b^0(x) \, \text{ a.e. in } \Omega \Bigg \},
\end{eqnarray*}
and the cost functions
\[
\theta_\nu(w^\nu,w^{-\nu}) := \frac{1}{2} \left \| \chi^\nu \left (\sum_{\mu=1}^{N} w^\mu \right ) + \chi^\nu  L^{-1} f- g^\nu \right \|^2_{L^2(\Omega)}
+ \frac{\gamma_\nu}{2 \beta_\nu^2} \left \| L w^\nu \right \|^2_{L^2(\Omega)},
\]
our elliptic multiobjective optimal control problem in the novel formulation $(III)$ 
writes as the jointly convex GNEP:
\[
\min_{w^\nu} \theta_\nu(w^\nu,w^{-\nu}) \quad \text{s.t.} \quad (w^\nu,w^{-\nu}) \in \widetilde{W} 
\]
for all $\nu=1,\ldots,N$. Solving this jointly convex GNEP gives us $(w^1,\ldots w^N)$ from which we can easily compute the state variable $y$ via \eqref{eq:y} and the controls $(u^1,\ldots,u^N)$ via \eqref{eq:u}, thus gaining the complete solution of our original problem (I). It was demonstrated in \cite{ADJG-JOTA}, that one can solve this reformulation (III) using a relaxation method that computes a best-response function and performs a line search exploiting a merit function, again based on the Nikaido-Isoda function.

\section{Lagrange multipliers, convex duality theory, and mixed formulations of nonsmooth variational problems}

In this section we provide mixed formulations of some nonsmooth variational problems and of associated variational inequalities. To achieve this goal we pursue a direct relatively simple approach to Lagrange multipliers that, however, heavily hinges on elliptic regularity theory. To put this approach in perspective we first shortly review the standard approach to Lagrange multipliers in convex duality theory that is based on the 
Hahn-Banach separation theorem.

\subsection*{A short review of convex infinite dimensional duality theory in function spaces}

The standard approach to prove existence of Lagrange multipliers for inequality constrained
optimization problems in infinite dimensional spaces is based on the 
Hahn-Banach separation theorem and thus needs interior point conditions, in particular a nonvoid interior of the ordering cone associated to the inequality constraint. 
In function spaces of continuous functions endowed with the maximum norm
with applications, e.g., to Chebychev approximation one can work with the topological interior of the ordering cone, see e.g., \cite{MR1429391}. However, the cone of 
non-negative $L^p$ functions and hence the ordering cone in the Sobolev spaces - relevant for pde constrained optimization - have empty topological interior. To overcome this difficulty one can resort to the concept of the so-called quasi-relative interior of a convex set introduced by Borwein and Lewis \cite{MR1167406}.
Therefore next we give the definition of this concept and  a short review of corresponding recent results on Lagrangean duality. 

Let $C$ be a nonvoid subset of a real normed space $X$. 
Let cl $C$, co $C$, cone $C$ denote the topological closure, 
convex hull, conical hull of $C$, respectively. 
Then for a given point $x \in C$,
the set 
$$
T_C(x) = \{y \in X: y =\lim_{n \to \infty} t_n(x_n - x), t_n > 0, x_n \in C 
~(\forall n \in \N), \lim_{n \to \infty} x_n = x  \}
$$
is called the {\it tangent cone (contingent cone)} to $C$ at $x$. If $C$ is convex, then 
$T_C(x) = \mbox{cl cone} (C- x)$.
With the dual space $X^\ast$ and the duality form $(.,.)$, the {\it normal cone} to $C$ at
$x \in C$ is defined by 
$$
N_C(x) = \{x^\ast \in X^\ast: (x^\ast,y-x) \leq 0, \forall y \in C  \} \,.
$$
Now the {\it quasi-interior} of a convex subset $C$ of $X$ is the set
$$ 
\mbox{ qi } C = \{x \in C: \mbox{cl cone }(C-x) =X \}
$$
and there holds the characterization, see \cite{MR2317768}, for 
$x$ in the convex set $C$:
$$
 x \in \mbox{qi } C \Leftrightarrow N_C(x)= \{0_{X^\ast}\} \,.
$$
Due to Borwein and Lewis \cite{MR1167406} is the following refinement 
of the notion of the quasi-interior: The {\it quasirelative interior}
 of a convex subset $C$ of $X$ is the set
$$ 
\mbox{ qri } C = \{x \in C: \mbox{cl cone }(C-x) 
\mbox{ is a linear subspace of } X \}
$$
and there holds the characterization, see \cite{MR2317768}, for 
$x$ in the convex set $C$:
$$
T x \in \mbox{qri } C \Leftrightarrow N_C(x)
\mbox{ is a linear subspace of } X^\ast  \,.
$$

These are useful concepts in $L^p$ function spaces 
and thus in Sobolev spaces as shown by the following example.

\begin{example}
Consider the Banach space $X=L^2 (T, \mu)$ with $1 \leq p < \infty$ 
on a measure space $(T,\mu)$ and the closed convex cone 
$C= \{ z \in X: z(t) \geq 0 \; \mu \mbox{ - a.e.} \; \}$. 
Then the characteristic function of $T$, $1= 1_T$ 
lies in qi $C$, hence in qri $C$. Indeed, by Lebesgue's theorem
of majorized convergence, any $x \in X$ can be approximated
by the sequence $\{x_n\}$ of truncations,
\[ x_n (t) = \left\{ \begin{array}{ll}
x(t) & \mbox{if $x(t) \geq - n$ a.e.}; \\
- n & \mbox{elsewhere}, 
\end{array} \right. \]
and clearly $x_n \in n (C-1)$.
\end{example} 

Now let us turn to inequality constrained convex optimization and 
Lagrangean duality theory. Consider the following primal optimization
problem:
$$
\mbox{(P)} \quad \inf_{x \in R} f(x) \,,
$$
where 
$$
R = \{x \in S: g(x) \in -C \},
$$
is assumed to be nonempty and $S$ a nonempty subset of $X$; $Y$ is another normed space partially ordered by a convex cone $C$; $f:S \to \R$ and $g:S \to Y$ 
are two maps such that the map $(f,g): S \to \R \times Y$, defined by 
$(f,g)(x) = (f(x),g(x)), \forall x \in S$   is convex-like with respect to the cone $\R_+ \times C \subset \R \times Y$, that is the set
$(f,g)S + \R_+ \times C$ is convex.  
Then the Lagrangian is
$$
L(x,\ell) = f(x) + (\ell,g(x)), \quad x \in S, \ell \in C^\ast  
$$
and the Lagrange dual problem to (P) reads
$$
\mbox{(D)} \quad \sup_{\ell \in C^\ast} \inf_{x \in S}
          [f(x) + (\ell,g(x))] \,,
$$
where $C^\ast = \{ x^\ast \in X^\ast : (x^\ast,x) \geq 0 , \forall x \in C \}$
is the dual cone to $C$.
While for the optimal values of (P) and (D), 
$\inf \mbox{(P)} = \inf_x \sup_\ell L(x,\ell)  \geq
\sup_\ell \inf_x L(x,\ell) =\sup \mbox{(D)}$ trivially holds,
one is interested in the equality of these optimal values 
and moreover in the existence of a {\it Lagrange multiplier},
that is, an optimal solution $\ell$ in (D).
This is called {\it strong duality}.

In the favorable situation when the topological interior of 
the ordering cone, $\mbox{int } C$, is not empty, 
the approach to strong duality in infinite dimensions
via the Hahn-Banach separation theorem
requires the easily verifiable Slater condition as 
a constraint qualification 
(see the important paper of Jeyakumar and Wolkowicz \cite{MR1167408}),
that is, the existence of a feasible point $\tilde x \in R$ such that 
$g(\tilde x) \in - \mbox{int } C$.

Thus one may be inclined to transfer this approach to the situation 
when the topological interior of $C$ is empty 
by replacing ``int'' by ``qri''.
However, this fails, as the following example due to
Daniele and Giuffr{\`e} \cite{MR2340691} shows.
   
\begin{example}
Let $X=S=Y=l^2, $ the Hilbert space of all real sequences $x = (x_n)_{n\in \N}$
with $\sum_{n=1}^\infty x_n^2 < \infty$ 
and $C = l^2_+$ the cone of all non-negative sequences in $l^2$.
Define $f: l^2 \to \R$ and $g: l^2 \to l^2$ respectively by
$$
f(x) =  \sum_{n=1}^\infty \frac{x_n}{n}; \quad  
(g(x))_n = - \frac{x_n}{2^n}, \quad \forall n \in \N \,.
$$
Then the feasible set $T= \{x \in l^2~|~-g(x) \in l^2_+ \} = l^2_+ $.
One has cl$(l^2_+ - l^2_+) = l^2, l^2_{++} := \mbox{ qri } l^2_+ = 
\{x \in l^2~:~x_n > 0, \forall n \in \N \}  \not= \emptyset$.
Take $\tilde x \in l^2_+, \tilde x_n = \frac{1}{n}$, then
$- (g(\tilde x))_n = \frac{1}{n 2^n}, -g(\tilde x) \in l^2_{++}$.
Further $\inf \mbox{(P)} = 0$ and $x = 0_{l^2}$ is the optimal solution of 
(P). On the other hand, for $\ell \in l^2_+$ we have
\begin{eqnarray*} 
\inf_{x \in l^2} [f(x) + (\ell,g(x))] 
&=&  \inf_{x \in l^2} \left [ \sum_{n=1}^\infty \frac{x_n}{n} -
             \sum_{n=1}^\infty \ell_n \frac{x_n}{2^n} \right ] \\
&=&  \inf_{x \in l^2} \sum_{n=1}^\infty \left [\frac{1}{n} - \frac{\ell_n}{2^n} \right ] x_n \\
&=& 
\left\{ \begin{array}{cc}
0 & \mbox{if } \ell_n = \ds \frac{2^n}{n} \, \forall n \in  \N  \,,\\[0.5ex] 
- \infty & \mbox{otherwise } \,.
\end{array} \right.
\end{eqnarray*}  
However, $\ell$ with $\ell_n = \frac{2^n}{n}$ does not belong to $l^2$.
hence $\sup \mbox{(D)} = - \infty$ and the optimal values do not coincide.
\end{example}

This example can also be given in a function space using the well-known isometry 
of $l^2$ and $L^2(0, 2\pi)$ based on Fourier expansion.

So in addition to a qri Slater-like condition one needs extra conditions to ensure strong duality. To this aim   
Bo{\c{t}},  Csetnek, and Moldovan \cite{MR2438595} 
introduce the following conic extension of (P) in the image space:
\begin{eqnarray*}
{\cal E}_{\inf\mbox{(P)}} &=& \{ (\inf\mbox{(P)} - f(x) - r, -g(x) -y):
x \in S, r \geq 0, y \in C \} \\
&=& ( \inf\mbox{(P)}, 0_Y ) - (f,g)~ S - \R_+ \times C \,,
\end{eqnarray*}
where as in classic convex duality theory only $\inf\mbox{(P)} \in \R$
is required, but not the existence of an optimal solution to (P).
Note that by feasibility of (P), $R \not= \emptyset$ implies 
$\inf\mbox{(P)} < \infty$
and in the case  $\inf\mbox{(P)} = - \infty$ strong duality    
trivially holds.

In this way Bo{\c{t}},  Csetnek, and Moldovan \cite{MR2438595} 
could prove the following strong duality result.

\begin{theorem} \cite[Theorem 4.1]{MR2438595}
Suppose that cl$(C-C) = Y$ and there exists some $\tilde x \in S$
such that $g(\tilde x) \in - \mbox{qri } C$. If  
\begin{equation} \label{BCM-CQ}
(0,0_Y) \notin \mbox{ qri co }
[{\cal E}_{\inf\mbox{(P)}} \cup \{(0,0_Y)\}] \,,
\end{equation} 
then strong duality holds.
\end{theorem}

\subsection*{A direct approach to Lagrange multipliers and dual mixed formulations of inequality constrained optimization and of VIs of the first kind}

We start with convex quadratic optimization in infinite dimensional spaces.
 Let $V$ be a real Hilbert space 
and let $Q$ be another real Hilbert space (for simplicity identified with its dual $Q'$).
Let $A \in {\cal L}(V,V')$ with $A = A', A \geq 0~$ 
(i.e., $\la Av,v\ra \geq 0,\forall v \in V$).
Further let $B \in {\cal L}(V,Q)$ and let $f \in V', g \in Q$ be fixed elements.
Moreover  let an order $\leq$ defined in $Q$ via a convex closed cone $C \subset Q$
 via $q \geq 0$, iff $q \in C$. 
With these data consider the convex quadratic optimization problem

\[
 (CQP) \qquad
    \left\{\begin{array}{l} \mbox{minimize }  f(v)=\frac{1}{2} \la Av,v \ra - \la f,v \ra
               \\ \mbox{subject to } Bv \leq g \:.
           \end{array}  \right.   
\]
					
This gives rise to the bilinear form $a(u,v):= \la Au,v \ra$
and the convex closed set 
$$ K(g) : = \{ v \in V~|~ Bv \leq g \}  \,,$$
which is translated from the cone
$$ K_0 : = \{ v \in V~|~ Bv \leq 0 \}.$$

As is well-known, a solution $u$ of $(CQP)$ is characterized
by the following  VI of the first kind - following the terminology in \cite{MR2423313}:
 
\[ (VI-1) \qquad
 u \in K(g),~ a(u,v-u) \geq \la f,v-u\ra,\, \forall v \in K(g) \,. 
\]

Here we present a simple approach - different from the approach reviewed above - to Lagrange multipliers. 
Assume that there exists a preimage of $g$ under $B$, $B \tilde g =g$. 
This allows to work with the duality on $V \times V'$,
obtain readily the existence of a Lagrange multiplier in the dual cone 
$$ K_0^+ = \{ \kappa \in V'~:~ \la \kappa, w\ra \geq 0, \, \forall w \in K_0 \} $$  
and arrive at the following characterization.

\begin{proposition} Let $u \in K(g)$. Then $u$ solves the above $(VI-1)$, iff there exists $\lambda \in K_0^+$ such that 
$(u,\lambda) \in V \times V'$ solves the mixed system

 \[(MP-1) \qquad
    \left\{\begin{array}{l} a(u,v) = \la \lambda,v \ra + \la f,v \ra     \\[0.5ex]
		\la \mu - \lambda, u-\tilde g \ra \geq 0  \:,
           \end{array}  \right.   \]

for all $v \in V, \mu \in K_0^+$. Further there holds the complementarity condition
$$ \la \lambda, u-\tilde g \ra = 0 \,. $$
\end{proposition}

\begin{proof}
Let $u \in K(g)$ solve the $(VI-1)$: 
$ a(u,v-u) \geq \la f,v-u\ra,\, \forall v \in K(g) $

Define $\lambda \in V'$ by
$\lambda(v) = a(u,v) -f(v)$. Then $(MP-1)_1$ holds. Further,
for any $v \in K_0$, $\tilde v := v + u$ lies in $K(g)$ and hence 
$$
 \lambda (v) =  a(u,\tilde v - u) - f(\tilde v - u) \geq 0 \,.
$$					
Thus $\lambda \in K_0^+$. Since $\tilde g \in K(g), u-\tilde g \in K_0$,
$$
\la \mu - \lambda, u - \tilde g \ra = \la \mu , u - \tilde g \ra -[a(u,u-\tilde g) -f(u-\tilde g)] \geq 0 $$
for any $\mu \in K_0^+$ and therefore $(MP-1)$ holds.

The complementarity condition follows from $(MP-1)_2$ by the choices
 $\mu = 2 \lambda$, $\mu = 0$.

Vice versa, let $v \in K(g)$, hence $v-\tilde g \in K_0$.
This implies by the complementarity condition
$$
\la \lambda, v - u \ra 
= \la \lambda, v - \tilde g \ra - \la \lambda, u - \tilde g \ra \geq 0 \,.
$$
Hence we arrive at
$$ a(u,v-u) = (f+ \lambda) (v-u) \geq f(v-u) \,. $$
\end{proof}
 
By the proof above it is clear that $(MP-1)$ is equivalent to the following complementarity
problem: Find $(u,\lambda) \in V \times V'$ such that 
\[(CP-1) \qquad
    \left\{\begin{array}{l} \lambda = A u - f    \\
		 \lambda \in K_0^+, u - \tilde g \in K_0 \\
		\la \lambda, u- \tilde g \ra = 0  \:.
           \end{array}  \right.   \]
					
Moreover the proof shows that the characterization above holds also with not necessarily symmetric bilinear forms, when the equivalence to convex quadratic optimization is lost; it even holds for nonlinear operators $A$ mapping a Banach space $V$ to its dual $V'$. 					
					
This approach applies to domain obstacle problems,  where the linear map $B$ is the imbedding map \cite{MR2424078}, say from $H^1(\Omega)$ to $L^2(\Omega)$ for linear scalar elliptic operators $L$ or more generally from $W^{m,p}(\Omega)$ to some $L^q(\Omega)$. 
It also applies to boundary obstacle problems or unilateral contact problems with the Signorini condition on some boundary part $\Gamma_c$ in appropriate function spaces, where the linear map $B$ is the trace map $\gamma$  \cite{MR1301021} to the boundary part $\Gamma_c$.
By this simple approach, the Lagrange multiplier lives in the dual of the Sobolev space of the variational problem, thus at first, is a general measure which may be singular.  Here regularity theory
- see the review in the second section of this paper -  comes into play to conclude that the Lagrange multiplier is indeed a $L^p$ function. Thus from an inequality constraint, one finally obtains a Lagrange multiplier $\lambda$ in the cone $L^p_+$ of non-negative $L^p$ functions on the domain $\Omega$. 
Thus we obtain the recent result \cite[Theorem 3.3]{MR3245967}
of Daniele, Giuffr{\`e}, Maugeri, and Raciti.  
When in the (scalar) mixed Signorini problem with a linear elliptic pde, there exists a multiplier $\ell$ to the inequality constraint
$\gamma  v \leq g \Leftrightarrow v|\Gamma_c \leq g$ a.e. 
that lives in the dual $Q'$ to the image space $Q = L^2(\Gamma_c)$, thus lies in  $L_+^2(\Gamma_c)$, then the multipliers $\ell$ and $\lambda$ are related by $ \lambda = \gamma^\ast \ell$, where $\gamma^\ast$ denotes the adjoint of the trace map $\gamma: H^1(\Omega) \to L^2(\Gamma_c)$  

Indeed, this direct simple approach to Lagrange multipliers and mixed formulations is used in an efficient numerical treatment of domain obstacle problems. Based on such mixed formulations the very effective biorthogonal basis functions with local support, due to Lamichhane and Wohlmuth \cite{MR2261019}, can be employed for approximation of the Lagrange multipliers in the hp-adaptive FEM for elliptic obstacle problems, see the recent paper \cite{Banz2011hpadaptive} of Banz and Schr\"{o}der.
 
\subsection*{A direct approach to Lagrange multipliers for VIs of second kind}

Here we consider non-smooth optimization problems of the form
$$ (NOP) \qquad \min \limits_{v\in V}  f(v)=\frac{1}{2} \la Av,v\ra - \la f,v \ra 
+ \varphi(v), 
$$
where $\varphi$ is convex, even positively homogeneous on $V$, but not differentiable in the classic sense. A prominent example encountered with given friction or Tresca friction
in solid mechanics  is

$$ \varphi_g (v) = \int_{\Gamma_c} g |v|~ds \quad (g \in L^\infty(\Gamma_c), g > 0)\,.$$

An optimal solution of $(NOP)$ is characterized as a solution of the VI of the second kind:
$$(VI-2) \qquad 
u\in V , \qquad \la Au,v-u\ra +  \varphi(v) - \varphi(u) \geq f(v-u),~\forall v \in V\,.$$

For the above example of $\varphi_g$ use  
\begin{eqnarray*}
\varphi_g ( v ) & = &  \int\limits_\Gamma g~| v | ~d \Gamma  =   \sup \left \{\int\limits_\Gamma g~ v~ \mu  ~~d \Gamma ~\Big |~ \mu \in L^2(\Gamma),
		|\mu| \leq 1 \right \},
\end{eqnarray*}

where $\sup$ is attained by $\mu = \mbox{sign } v$,
set $$
M := \{ \mu \in L^2(\Gamma), |\mu| \leq 1\}
$$
and arrive - as it is shown in more general terms in Proposition \ref{prop-2} 
below - at the mixed problem:

Find $u \in V = H^1 ( \Omega), \lambda \in M$ such that for all $v \in V, \mu \in M$

\begin{eqnarray*}
 \left\{    \begin{array} {lll}
 \la A u , v \ra +   
\int\limits_\Gamma g~ v~ \lambda  ~d \Gamma & = & \la f,v \ra  \,,
 \\[2ex]
\int\limits_\Gamma g~ u~ (\lambda - \mu)  ~d \Gamma & \geq & 0 \,.
\end{array} \right. 
\end{eqnarray*}

To reveal the duality structure, introduce
$$
M(g) = \{ \mu \in  L^2(\Gamma), |\mu| \leq g  \mbox{ a.e. }\} \,.
$$
Although, with $g \in L^\infty(\Gamma)$, this set is clearly contained in $L^\infty(\Gamma)$,
we stick to the easier treatable $L^2$ duality. Thus $(VI-2)$ is equivalent   
-  as it is shown in more general terms in Proposition \ref{prop-2} 
below - to the mixed problem:  

Find $u \in V = H^1 ( \Omega), \lambda \in M(g)$ 
such that for all $v \in V, \mu \in M(g)$

\begin{eqnarray*}
 \left\{    \begin{array} {lll}
 \la A u,v \ra_{V^\ast \times V}  +    \la \lambda, v \ra_{L^2(\Gamma)}  
& = & \la f,v \ra_{V^\ast \times V}  \,,
 \\[1ex]
\la u, \lambda - \mu \ra_{L^2(\Gamma)}  & \geq & 0 \,.
\end{array} \right. 
\end{eqnarray*}

Indeed, in the more general setting of a reflexive Banach space $V$ , a map 
$A: V \to V^\ast, f\in V^\ast$ and a sublinear functional $\varphi: V \to \R$,
we have the following result using the convex weakly-$\ast$-compact subdifferential
$$
P:= \partial \varphi(0) = 
\{ q \in V^\ast,  \la q,v \ra \leq \varphi(v), \forall v \in V \} \,.
$$

\begin{proposition} \label{prop-2}
 $u \in V$ solves the above $(VI-2)$, iff there exists
 $p \in P$ such that 
$(u,p) \in V \times V^\ast$ solves the mixed system

 \[(MP-2) \qquad
    \left\{\begin{array}{l} \la A u,v \ra +  \la p,v \ra = \la f,v \ra     \\[0.5ex]
		\la p - q, u \ra \geq 0  \:,
           \end{array}  \right.   \]

for all $v \in V, q \in P$.
 \end{proposition}

\begin{proof}
Let $u \in V$ solve  $(VI-2)$. Then the choice $v=0$ gives 
\begin{equation} \label{prop_21}
\la Au,u\ra +  \varphi(u) \leq f(u) \,,
\end{equation} 
whereas  the choice $v=t w,  w \in V, t > 0, t \to \infty$ gives 
for all $w \in V$, 
\begin{equation} \label{prop_22}
\la Au,w\ra +  \varphi(w) \geq f(w) \,,
\end{equation} 
hence from (\ref{prop_21}) and (\ref{prop_22}) we get
\begin{equation} \label{prop_23}
\la Au,u\ra +  \varphi(u) = f(u) \,.
\end{equation} 
Note that (\ref{prop_23}) and (\ref{prop_22}) imply $(VI-2)$, 
hence these assertions are equivalent to $(VI-2)$.

Define $p \in V^\ast$ by
$p  = f - A u$. Then $(MP-2)_1$ trivially holds. Further from (\ref{prop_22}),
for any $w \in V$, $\varphi(v) \geq \la p,w \ra$, hence
$p$ lies in $\partial \varphi(0) =P$. Finally from  (\ref{prop_23}),
$ \varphi(u) = \la p,u \ra$, hence  $(MP-2)_2$ follows.

Vice versa,  $(MP-2)_2$ implies $ \varphi(u) = \la p,u \ra$,
hence together with $(MP-2)_1$ and the choice $v=u$ gives  (\ref{prop_23}). 
Since $\varphi(v) \geq \la p,v \ra$, from $(MP-2)_1$ we arrive at (\ref{prop_22}).
\end{proof}

Similarly as discussed in the previous subsection, the regularity 
of the multiplier $p$ hinges on the regularity of the datum $f$ and in particular on the regularity of the solution $u$ of the $(VI-2)$
via the map $A$.
 
To apply the above general result to the friction-type functional $\varphi_g$ 
we only have to set $V = H^1(\Omega), \varphi = \varphi_g \circ \gamma $
with the linear continuous trace operator $\gamma$ that maps
$H^1(\Omega)$ onto $H^{\frac{1}{2}}(\Gamma)$ dense in $L^2(\Gamma)$
and use the subdifferential chain rule \cite{MR0373611}
$$
\partial \varphi = 
\partial (\varphi_g \circ \gamma) = \gamma^\ast \partial \varphi_g \gamma \,.
$$
Note that this chain rule holds as an equality, since $\varphi_g$ is real-valued and so the constraint qualification $0 \in \mbox{int }(\mbox{range } \gamma - \mbox{dom } \varphi)$
is trivially satisfied.

To conclude this subsection let us mention other duality relations and mixed formulations useful in numerical treatment of variational inequalities of the second kind.
By $(L^1, L^\infty)$ duality and density one obtains 
$$ \varphi_g (v) = \int_{\Gamma_c} g |v|~ds  = \sup \left \{ \int_{\Gamma_c} g~ v~ \mu~ds
~ \Big |~ {\bf \mu \in C(\Gamma)}, |\mu| \leq 1 \right \}  \,.$$
This is used in convergence proof of Finite Element Methods and Boundary  Element Methods, see 
\cite{MR2566769,MR3061074}. 

Another way to cope with the  nondifferentiable functional $\varphi_g$
is to decompose the modulus function $|\rho| = \rho^+ + \rho^-$ with the positive part 
$\rho^+  = \max( \rho,0) \geq 0$ 
and the negative part $\rho^- = \max( - \rho,0) \geq 0$. This leads to inequality constrained problems considered in the previous subsection
what is not elaborated here further.

\subsection*{A direct approach to Lagrange multipliers for more general VIs}

To conclude this section we deal with the more general   
$$
(VI-3) \qquad 
u\in K , \qquad \la A(u),v-u\ra +  \varphi(v) - \varphi(u) \geq f(v-u),~\forall v \in K\,,
$$
where as above $f \in V^\ast$, $\varphi: V \to \R$ is sublinear and now
$K \subset V$ is a convex closed cone with vertex at zero. 
A VI of this form occurs in unilateral contact of a linear elastic body with a rigid foundation under the Tresca friction law, if the initial gap between body and foundation is zero, see \cite{MR3061074}. 
The more general setting, in particular for non-zero gap, with $K$ and $\varphi$ convex would encompass also the VIs of first kind studied before, but needs additional
arguments. Therefore we prefer this simpler homogeneous setting 
to elucidate the direct approach to Lagrange multipliers.
In this setting we have the following result in a general 
locally convex topological vector space $V$ for a not necessarily linear operator
$A: V \to V^\ast$. 

\begin{theorem} \label{prop-3} Let $u \in K$. Then
 $u$  solves the above $(VI-3)$, iff there exist
 $p \in P= \partial \varphi(0)$ and $\lambda \in K^-$ such that
the complementarity condition $\la \lambda, u \ra = 0$ holds and 
$(u,p,\lambda) \in V \times V^\ast \times V^\ast$ solves the mixed system
 \[(MP-3) \qquad
    \left\{\begin{array}{l} \la A (u),v \ra +  
\la p + \lambda ,v \ra = \la f,v \ra     \\[0.5ex]
\la p - q, u \ra \geq 0  \:,
           \end{array}  \right.   \]
for all $v \in V, q \in P$.
 \end{theorem}

\begin{proof}
Let $u \in K$ solve  $(VI-3)$. Then 
we first proceed as in the proof of Proposition \ref{prop-2}.
The choice $v=0$ gives 
\begin{equation} \label{prop_31}
\la A(u),u\ra +  \varphi(u) \leq f(u) \,,
\end{equation} 
whereas  the choice $v=t w,  w \in K, t > 0, t \to \infty$ gives 
for all $w \in K$, 
\begin{equation} \label{prop_32}
\la A(u),w\ra +  \varphi(w) \geq f(w) \,,
\end{equation} 
hence from (\ref{prop_31}) and (\ref{prop_32}) we get
\begin{equation} \label{prop_33}
\la A(u),u\ra +  \varphi(u) = f(u) \,.
\end{equation} 
Note that (\ref{prop_33}) and (\ref{prop_32}) imply $(VI-3)$, 
hence these assertions are equivalent to $(VI-3)$.

Define $\ell \in V^\ast$ by
$\ell  = f - A(u)$. From (\ref{prop_32}) we find
\begin{equation} \label{prop_34}
\varphi(v) \geq \la \ell,w \ra\,, \quad  \forall w \in K \,. 
\end{equation}
Now we claim that $\ell \in K^- + P$ and that hence $(MP-3)_1$ holds. Note that both $K^-$ and $P$ are convex closed sets, moreover $P$ is weakly* compact in $V^\ast$.
So the claim can be shown by an indirect argument employing the separation theorem. Here we use that 
$$ \varphi(w) = \max_{q \in P} \la q,w \ra $$
and thus (\ref{prop_34}) means that for any $w \in K$ there exists 
$q \in P$ such that $\la q,w \ra \geq \la \ell,w \ra$.
Therefore by the extension lemma \cite[Theorem 2.2]{MR821398} 
(which is a refined version of the famous Fan-Glicksberg-Hoffman theorem of alternative and is proved from a fixed point theorem or from the separation theorem)
there exists $p \in P$ such that  $\la q,w \ra \geq \la \ell,w \ra$ holds for all $w  \in K$. Now define $\lambda = \ell - p$, hence $\lambda \in K^-$
and $\ell = \lambda + p \in K^- + P $ as claimed.

From (\ref{prop_33}) we obtain
$$
\la p,u \ra  \leq \varphi(u) = \la \lambda + p, u \ra \,,
$$  
hence $\la \lambda,u \ra \geq 0$. Since $u \in K, \lambda \in K^-$, 
the complementarity condition $\la \lambda,u \ra = 0$ follows.
This gives with (\ref{prop_33}) that
$ \varphi(u) = \la p,u \ra$, hence  we arrive at $(MP-3)_2$.

Vice versa,  $(MP-3)_1$ implies together with the complementarity condition  and $(MP-3)_2$  
$$ \la f - A(u), u \ra = \la \lambda + p ,u \ra = \la p,u \ra = \varphi(u)\,,
$$
hence  (\ref{prop_33}). In view of $\lambda \in K^-$ we conclude from
 $(MP-3)_1$  that for any $w \in K$,
$$ \la f - A(u), w \ra = \la \lambda + p ,w \ra 
\leq  \la p,w \ra  \leq \varphi(w)\,,
$$  
hence (\ref{prop_32}).
\end{proof}

\section{Conclusions and outlook}

We have seen the crucial role of elliptic regularity theory in two instances.
First with  elliptic multiobjective optimal control formulated as jointly convex GNEP the regularity of the solution of the underlying pde was needed to arrive at a reformulation that was the basis for an efficient numerical solution method.
In this approach we had to require that the domain where the elliptic pde lives is convex or sufficiently smooth. On the other hand, real-world domains may have reentrant corners or are only piecewise smooth. This leads to the question how this approach can be refined using the well-known elliptic theory in nonsmooth domains \cite{MR3396210},
abandoning classic Sobolev spaces, and working instead with weighted Sobolev spaces
\cite{MR926688}.

Then we have presented a direct approach to Lagrange multipliers in inequality constrained and related nonsmooth boundary value problems which gives an immediate link between the regularity of the Lagrange multiplier and the regularity of the solution of the problem.
As already the one-dimensional obstacle problem demonstrates, 
there is a threshold of smoothness, however, that in general cannot be overstepped even if the data are arbitrarily smooth. The regularity theory for frictionless unilateral contact reported from the work \cite{MR656002} has shown the influence of the switching points, where the boundary conditions change, on the smoothness of the solution. So one may be interested in a more detailed analysis in weighted Sobolev spaces \cite{MR926688} that takes the switching points in account.

Finally let us point out that we have here considered frictionless monotone unilateral contact problems. Nonmonotone contact problems can be put in primal form as hemivariational inequalities (HVIs).
While the theory of HVIs is well developed, the numerical solution of these problems is in its infancy; here we can refer to 
\cite{MR1365370,MR2817476,MR3113542,MR3245969,MR3335209,MR3355947}
(ordered according to publication date).
So one may ask for mixed formulations with appropriate Lagrange multipliers that would allow the development of mixed finite element procedures for the efficient solution of these nonconvex variational problems.

\bibliographystyle{abbrv}
\bibliography{BIB_RegNumSol}
 
\end{document}